\documentclass[conference]{IEEEtran}
\ifCLASSINFOpdf
  \usepackage[pdftex]{graphicx}
\else
\fi
\usepackage{amsmath}
\usepackage{amsthm}
\usepackage{array}

\ifCLASSOPTIONcompsoc
  \usepackage[caption=false,font=normalsize,labelfont=sf,textfont=sf]{subfig}
\else
  \usepackage[caption=false,font=footnotesize]{subfig}
\fi
\usepackage{url}

\newcommand{\expect}[1]{\mathrm{E}\left[#1\right]}

\begin{document}
%
\title{Optimal Orderings of $k$-subsets for Star Identification}

\author{\IEEEauthorblockN{Joerg H. Mueller\IEEEauthorrefmark{1}, Carlos S\'{a}nchez-S\'{a}nchez\IEEEauthorrefmark{1},
Lu{\'i}s F. Sim{\~o}es\IEEEauthorrefmark{2} and
Dario Izzo\IEEEauthorrefmark{1}}
\IEEEauthorblockA{\IEEEauthorrefmark{1}Advanced Concepts Team,
European Space Agency,
Noordwijk, The Netherlands\\ Email: \{joerg.mueller, carlos.sanchez, dario.izzo\}@esa.int}
\IEEEauthorblockA{\IEEEauthorrefmark{2}Computational Intelligence Group,
Vrije Universiteit Amsterdam, Amsterdam, The Netherlands\\
Email: luis.simoes@vu.nl}}


\maketitle

\begin{abstract}
Finding the optimal ordering of $k$-subsets with respect to an objective function is known to be an extremely challenging problem. In this paper we introduce a new objective for this task, rooted in the problem of star identification on spacecrafts: subsets of detected spikes are to be generated in an ordering that minimizes time to detection of a valid star constellation.
We carry out an extensive analysis of the combinatorial optimization problem, and propose multiple algorithmic solutions, offering different quality-complexity trade-offs.
Three main approaches are investigated: exhaustive search (branch and prune), goal-driven (greedy scene elimination, minimally intersecting subsets), and stateless algorithms which implicitly seek to satisfy the problem's goals (pattern shifting, base unrank).
In practical terms, these last algorithms are found to provide satisfactory approximations to the ideal performance levels, at small computational costs.
\end{abstract}


%
\IEEEpeerreviewmaketitle

\section{Introduction}

In this paper, we introduce a novel objective for the generation of all $k$-subsets of $n$ elements and we discuss the structure of the resulting combinatorial optimization task. In general, improving the order of elements in a sequence towards some objective is recognized to be a complex optimization task \cite{dewar} with interesting applications in computer science such as unit test coverage \cite{bryce,rothermel}.

\begin{figure}[tbp]
    \centering
    \includegraphics[width=0.8\columnwidth]{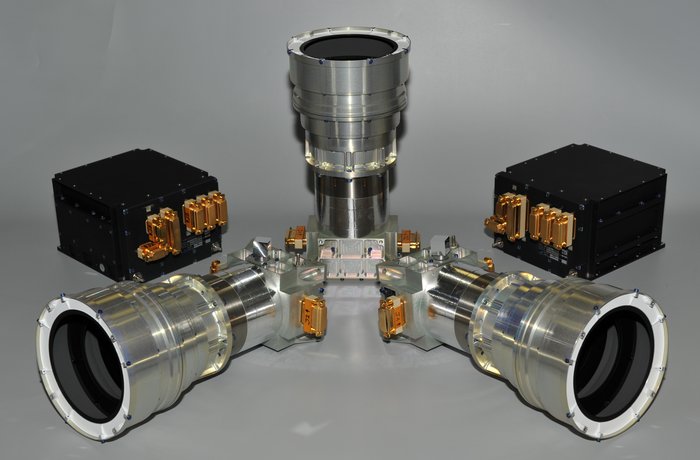}
    \caption{The flight model of the hydra star-tracker currently flying on the ESA Sentinel 3A satellite mission.}
    \label{fig:hydra}
\end{figure}

Star trackers (see Figure \ref{fig:hydra}) are a common sensor used by a spacecraft to determine its attitude by looking at fixed stars. The problem we here formalize and tackle was suggested by the work of Mortari et al. \cite{MortariEtAl2004} on the design of efficient algorithms for star identification. In that paper, the authors consider a ``lost-in-space" spacecraft attitude identification problem: find the orientation of a spacecraft in deep space using a single star tracker image. Such a problem corresponds to that of identifying $k$ stars in a ``scene" (i.e. a picture taken by the star tracker) containing $n$ spikes of which $t$ are unknown stars and the rest are artifacts due to various disturbances present in harsh space environments. Part of the algorithm proposed in that paper, called the Pyramid Algorithm and today widely used in many star-trackers in orbit, needs to generate all $\binom{n}{3}$ combinations in a smart order that allows the discovery of three true stars in the scene from a minimal number of star catalog queries. 

\section{Background}

A star tracker utilizes a star catalog containing the positions of all known stars having brightness larger than some predefined threshold. A camera is used as sensory device to capture the light coming from the stars which will then create spikes in the image. The crucial software part of a star tracker is the star identification algorithm whose task is to match spikes extracted from the image with stars in the star catalog. With the matched stars it is then possible to calculate the orientation of the spacecraft. One major problem for star identification algorithms is that not all spikes in an image are associated to real stars. Spikes can be caused by reflections on debris, radiation or other spurious sources. A correct and fast star identification is critical in space: a failure of the spacecraft to detect its own attitude promptly may lead to the complete failure of the mission.

\subsection{The Pyramid Algorithm}

A widely used star identification algorithm \cite{SpratlingEtAl2009} is the Pyramid Algorithm developed by Mortari et al. \cite{MortariEtAl2004}. The algorithm is based on the distances between two stars in an image. A database stores the distance between each star pair from the catalog up to the maximum distance that is possible within the field of view of the camera. The look-up of an arbitrary distance within the measurement tolerance typically results in a number of possible star pairs from the catalog. It is therefore necessary to build higher order graphs with the spikes in the image, such as a triangle consisting of three distances. The resulting star pairs from the distance look-ups in the database have to be matched to form an actual triangle given the IDs of the stars in the star catalog. 

\begin{figure}[tbp]
    \centering
    \includegraphics[width=0.8\columnwidth]{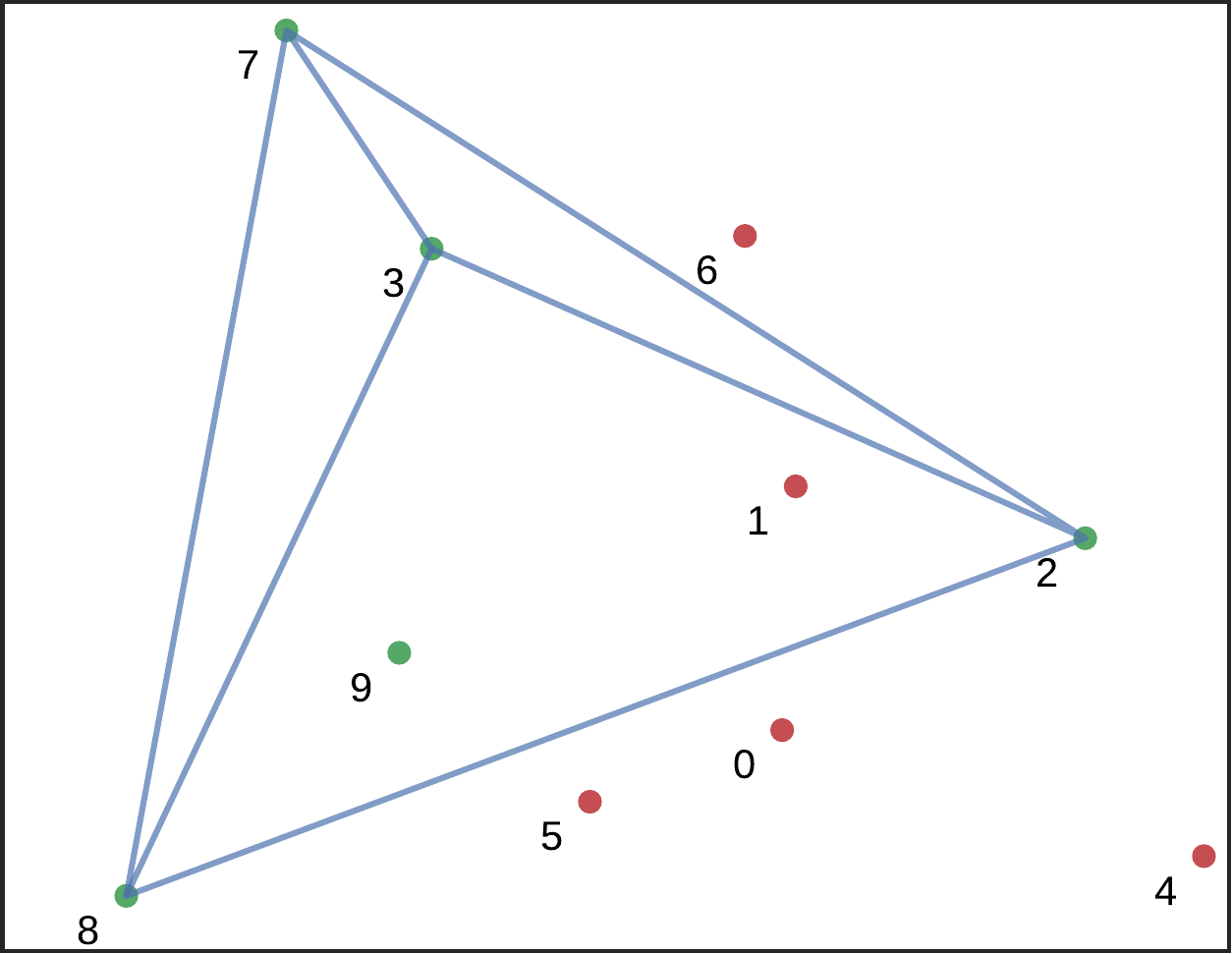}
    \caption{A random scene as a star tracker would see it. Ten ($n=10$) detected spikes in the image are numbered arbitrarily from 0 to 9 and marked as circles. Artifacts are colored in red, actual stars ($t=5$) are colored in green. The pyramid algorithm detected the four stars connected with blue lines as a valid pyramid.}
    \label{fig:scene1}
\end{figure}

An important part of the algorithm is to select spike triplets resulting in a sequence of queries to the database able to find a positive match as quickly as possible. This problem equals to that of wisely selecting the order of all possible $k$-subsets of $n$ given spikes, where, in this case, $k=3$.

For example, as discussed also in the original paper by Mortari, the lexicographic order (which tries the spikes in the sequence 0-1-2, 0-1-3, 0-1-4, etc.) is clearly problematic, as the first spike is repeated until all possible combinations containing it have been also tried. In the case of such a spike being spurious this approach results in many unnecessary database queries. As part of the Pyramid Algorithm, Mortari et al. proposed an algorithm to generate the subsets in a smarter fashion, while acknowledging that the mathematically optimal solution is of interest, but unknown. This stateless algorithm, which we call \emph{pattern shifting} algorithm, is described in the original paper and here reported for convenience:

\begin{verbatim}
for dy from 1 to n-2
    for dz from 1 to n-dy-1
        for x from 0 to n-dy-dz-1
            y = x + dy
            z = y + dz
            next combination is [x, y, z]
\end{verbatim}

This code would produce the sequence 0-1-2, 1-2-3, 2-3-4, etc.

\section{Problem Description}

Consider a star tracker camera retrieving an image containing $n$ spikes arbitrarily numbered from 0 to $n-1$ of which an unknown subset $s = \{p_1, p_2, \dots, p_t\}$ of cardinality $t$ are actual stars and the remaining spikes artifacts. A fixed $n$, $s$ defines what will here be called a scene, which can be thought of as the underlying, unknown, ground truth. A query $q = \{c_1, c_2, \dots, c_k\}$ is a $k$-subset of the spikes and is said to discover the scene if $q \subseteq s$.

We study the problem of generating, for a fixed $n$, an ordered sequence $Q=\{q_j, j\in 1..M\}$ of $M$ queries able to discover the unknown scene, on average, at the smallest trial $j$. To guarantee that $j$ always exists, we only consider a sequence $Q$ if it contains all $M = N = \binom{n}{k}$ possible queries (in which case we refer to $Q$ as \emph{complete}) and we assume $t\ge k$ as scenes with less than $t$ spikes cannot be discovered. The set $S = \{s_1, s_2, \dots \}$ for a fixed $n$ contains all possible scenes that are discoverable ($t \ge k$). Due to lack of prior knowledge about $S$, we assume all scenes to have equal probability $|S|^{-1}$, where:
\begin{equation}
    |S| = \sum_{t=k}^{n} \binom{n}{t}
\end{equation}
We formalize the problem as:

\begin{equation}\label{eq:ttd}
\begin{array}{rl}
\mbox{find:} & Q \\
\mbox{to minimize:} &  T(Q) = \frac{1}{|S|}\sum_{s \in S} \tau\left(Q, s\right)
\end{array}
\end{equation}
where $T(Q)$ is the expected time to discovery, that is the average across all possible scenes $s \in S$, of $\tau(Q, s)$: the index $j$ of the first $q_j \in Q$ discovering $s$.

\subsection{Discovery of Scenes}

Let $D(q_i)$ be the number of scenes discovered by the query $q_i$ but not by any of the previous queries $q_{j< i}$. Clearly, for each of such scenes, $\tau(Q, s)=i$. We may then rearrange the terms in the expression for the expected time to discovery in Eq.(\ref{eq:ttd}) grouping together all scenes discovered by $q_i$ and thus summing over all $q_i \in Q$:
\begin{equation}\label{eq:TQq}
T(Q) = \frac{1}{|S|} \sum_{i=1}^{|Q|} i D(q_i).
\end{equation}
In this form, the objective function suggests that the number of scenes discovered at each $i$ should be larger than those discovered at any $j>i$, a simple thought that will form the basis of some of the most successfull algorithms here introduced.

\subsection{Example}

Consider a simple scenario with $n=5$ spikes and $k=3$. We assume the spikes are numbered from 0 to 4. Let us compute the time to discovery $T$ of the sequence $Q$ generated by the pattern shifting algorithm by Mortari et al.:

\begin{equation*}
    Q =
    \begin{Bmatrix}
        \{ 0, 1, 2 \}, \\
        \{ 1, 2, 3 \}, \\
        \{ 2, 3, 4 \}, \\
        \{ 0, 1, 3 \}, \\
        \{ 1, 2, 4 \}, \\
        \{ 0, 1, 4 \}, \\
        \{ 0, 2, 3 \}, \\
        \{ 1, 3, 4 \}, \\
        \{ 0, 2, 4 \}, \\
        \{ 0, 3, 4 \} \\
    \end{Bmatrix}
\end{equation*}

As a start, consider the scene $s_e=\{1, 2, 3, 4\}$. As $q_1 = \{0,1,2\}\not \subseteq s_e$ does not discover the scene, but $q_2 =\{1,2,3\}\subseteq s_e$ does, the number of queries necessary to discover this particular scene is $\tau\left(Q, s_e\right) = 2$. To compute the total time to discovery $T(Q)$, we must sum over all possible scenes (of which there are 16 with $t \ge 3$ : one with $t=5$, five with $t=4$ and ten with $t=3$). Listing all possibilities it is not difficult to find that we have four scenes discovered by $q_1$, two discovered for each $q_2$ to $q_4$ and one for each $q_5$ to $q_{10}$. This, according to Eq.(\ref{eq:TQq}) gives a final score of:
$$
T(Q) = \frac{1}{16} \left( 4\cdot1 + \sum_{i=2}^{4} 2 \cdot i + \sum_{i=5}^{10} 1\cdot i \right) = \frac{67}{16}
$$
The optimum in this case could be reached, for example, by swapping $q_2$ and $q_{10}$ leading to  a score of $T(Q) = \frac{65}{16}$. The lexicographic sequence, for reference, has, in this case, a score of $T(Q) = \frac{71}{16}$.

\section{Analysis}\label{sec:analysis}

We discuss some properties useful to design algorithms aimed at solving the problem stated in Eq.(\ref{eq:ttd}).

\subsection{Equivalent sequences}\label{sec:uniqueness}

No sequence $Q$ has a unique score. Consider the equivalence class $[Q]$ containing all sequences having the same time to discovery. The cardinality of this set is at least $n!$. Since this property is also valid for the optimal sequence, we conclude that there are at least $n!$ solutions. This result follows immediately noting that any permutation of the $n$ spike IDs does not change the score. Therefore, applying the same permutation to all elements in each $q_i \in Q$ will result in an equally scored $Q'$.

\subsection{Scene Discovery Monotonicity}\label{sec:monotonicity}

If $Q$ is optimal, then the values $D(q_i)$ monotonically decrease with $i$. This property, called the monotonicity property, states that the number of scenes discovered by each successive query in an optimal sequence is a monotonically decreasing sequence.

It is illustrated in Table \ref{table:monotonicity}. The sequence on the left cannot be optimal as it violates the monotonicity property. By swapping the highlighted queries we obtain a better sequence.

\begin{table}
\centering
\caption{Monotonicity of the Optimal Order}
\label{table:monotonicity}
\def\arraystretch{1.1}
\begin{tabular}{cc|cc}
$q_i$ & $D(q_i)$ & $q^*_i$ & $D(q^*_i)$ \\ \hline
\{0, 1, 2\}                & 4                              & \{0, 1, 2\}                 & 4                           \\
\{0, 3, 4\}                & 3                              & \{0, 3, 4\}                 & 3                         \\
\textbf{\{0, 1, 3\}}                 & \textbf{1}                              & \textbf{\{2, 3, 4\}}                 & \textbf{2}                           \\
\textbf{\{2, 3, 4\}}               & \textbf{2}                              & \textbf{\{0, 1, 3\}}            & \textbf{1}                            \\
\{1, 2, 4\}                & 1                              & \{1, 2, 4\}                 & 1                            \\
\{0, 1, 4\}                & 1                              & \{0, 1, 4\}                 & 1                            \\
\{0, 2, 3\}                & 1                              & \{0, 2, 3\}                 & 1                            \\
\{1, 3, 4\}                & 1                              & \{1, 3, 4\}                 & 1                            \\
\{0, 2, 4\}                & 1                              & \{0, 2, 4\}                 & 1                            \\
\{1, 2, 3\}                & 1                              & \{1, 2, 3\}                 & 1    \\ \cline{2-2} \cline{4-4} 
& $T(Q)$ = 64 & & $T(Q^*)$ = 63
\end{tabular}
\end{table}

\begin{proof}

Consider two consecutive queries $q_i, q_{i+1} \in Q$ such that $D(q_{i}) < D(q_{i+1})$, we show that swapping the subsets $q_i, q_{i+1}$ results in a new sequence $Q^*$, where $q^*_{i} = q_{i+1}$ and $q^*_{i+1} = q_{i}$, with a lower average time to discovery implying that, for the optimal sequence, $D(q_i)$ monotonically decreases with $i$.

Due to the swap, $D(q_{i}) \ge D(q^*_{i+1})$ because of scenes that can be discovered by both $q_i$ and $q_{i+1}$. According to Eq. \eqref{eq:TQq}, The contribution of the queries $q_i$, $q_{i+1}$ to $T(Q)$ is:
$$T_i = i D(q_i) + (i+1) D(q_{i+1})$$
$T^*_i$ can be similarly defined as the contribution of the queries $q^*_i$, $q^*_{i+1}$ to $T(Q^*)$. Using $D(q_{i}) + D(q_{i+1}) = D(q^*_i) + D(q^*_{i+1})$, as both queries together always discover the same amount of scenes, we get
$$T_i - T^*_i = D(q_{i+1}) - D(q^*_{i+1}) > 0,$$
as $D(q_{i+1}) > D(q_{i}) \ge D(q^*_{i+1})$.

Swapping the two queries does not alter the contributions of the previous and following queries. Therefore, $T_i > T^*_i$ results in $T(Q^*) < T(Q)$.

\end{proof}

\subsection{The Average $T$}

Consider the set $\mathcal Q$ of all sequences that are \emph{complete} and do not contain repetitions. Following Eq.(\ref{eq:ttd}) and Eq.(\ref{eq:TQq}) we may restate the fundamental problem subject of this paper as:
\begin{equation}
\begin{array}{rl}
\mbox{find:} & Q \in \mathcal Q\\
\mbox{to minimize:} &  T(Q) = \frac{1}{|S|} \sum_{i=1}^{|Q|} i D(q_i).
\end{array}
\end{equation}
The dimension of our search space, that is the cardinality of $\mathcal Q$, can be computed noting that the number of queries $q_i \in Q$ is $N = \binom{n}{k}$. The total number of different sequences in $\mathcal Q$, the search space dimension, is then $N!$. It is of interest to compute $\sigma$, the average across the whole search space of $T(Q)$. Any algorithm producing a sequence scoring less than $\sigma$ will be considered as making good use of the problem structure, while the opposite can be said for algorithms that generate sequences scoring more than $\sigma$. The value $\sigma$ can be computed as follows:
\begin{equation*}
    \sigma = \frac{1}{N!} \sum_Q T(Q)
\end{equation*}
Computing $\sigma$ directly from this definition is quite expensive, thus we derive a simpler formula that has a much lower complexity. Let us start by assuming the first $i-1$ queries in some $Q$ have not discovered a scene $s_t$ having $t$ true stars. The probability that the following query $q_i \in Q$ will discover it is denoted with $p(q_i, t)$. This probability is  the fraction between the number of queries able to discover $s_t$ and the total number of remaining, possible queries:
\begin{equation*}
    p\left(q_i, t\right) = \frac{\binom{t}{k}}{N - i + 1},
\end{equation*}
This probability increases with $i$ and is well defined only when $i \le N+1-\binom{t}{k}$.
We may then compute the unconditional probability $p_i(t)$ that $s_t$ will be detected at the $i$-th query $q_i$ as:
\begin{equation*}
    p_i(t) = p\left(q_i, t\right) \prod_{j=1}^{i-1} \left(1 - p\left(q_j, t\right)\right)
\end{equation*}
Using this equation, it can be shown that $\sigma$ can be derived as:
\begin{equation}\label{eq:sigma}
   \sigma = \frac{1}{|S|} \sum_{t=k}^{n} \left( \left|S_t\right| \sum_{i=1}^{N} i p_i(t) \right),
\end{equation}
where $\left|S_t\right| = \binom{n}{t}$ is the number of possible scenes with exactly $t$ true stars.
The complexity of computing this is $O(N n)$. Some examples are shown in Table~\ref{tab:ET}, showing that the expectation for a fixed $k$ has a maximum at some $n$.

\begin{table}
\centering
\caption{Values of $\sigma$ for Different $n$ and $k$.}
\label{tab:ET}
\begin{tabular}{c|c|r}
$n$ & $k$ & \multicolumn{1}{c}{$\sigma$} \\ \hline
5 & 3 & 4.2 \\
10 & 3 & 17.4 \\
20 & 3 & 16.0 \\
50 & 3 & 9.9 \\
100 & 3 & 8.8 \\
\end{tabular}
\hspace{1em}
\begin{tabular}{c|c|r}
$n$ & $k$ & \multicolumn{1}{c}{$\sigma$} \\ \hline
20 & 1 & 2.0 \\
20 & 2 & 5.1 \\
20 & 5 & 322.5 \\
20 & 10 & 32528.1 \\
20 & 15 & 5748.7 \\
\end{tabular}
\end{table}

\section{Algorithms}

We propose several algorithms to generate a solution to the problem \eqref{eq:ttd} that we group in three families. First we detail algorithms able to find an optimal solution. Since these reveal to be computationally intractable with growing $n$, we introduce a second family including algorithms making use of different heuristics to search in $\mathcal Q$ resulting in a great reduction of the computational cost but without the guarantee of finding the optimal sequence. Finally, the algorithms of the last family are stateless $k$-subset (queries) generators. They have a considerably lower computational cost, making them attractive for real-time generation in star trackers. This family includes a generalization of the pattern shifting algorithm described by Mortari et al. \cite{MortariEtAl2004}. 

\subsection{Finding an optimal solution}

We start with the naive brute-force approach computing the score for all $Q \in \mathcal Q$ and returning the minimum. The same result, implemented as a branch and prune algorithm making use of the properties introduced in section \ref{sec:analysis} to prune large portion of the search space not containing the optimal solution.

\subsubsection{Brute-force}

In the brute-force approach one simply determine the time to discovery for all possible sequences of queries. This way of determining the best sequence comes at a high cost, with a time complexity of $O\left(N! k |S|\right)$.
Every permutation of the $N$ queries has to be checked against every possible scene. It is easy to understand that this complexity grows out of any reasonable bounds even for low numbers of $n$ and $k$.

\subsubsection{Branch and prune}

The brute-force algorithm can be implemented as branch and prune tree search. Each level in the tree adds another query to the sequence. This means the root node has $N$ neighbors and the branching factor at each level is decreasing by one from there, resulting in a depth of $N$ for all the leaves in the tree.

To optimize computational performance, scenes and queries can be stored and processed in a binary representation. The bit-wise \emph{AND} operation allows to check if a query would discover a scene. The search for the optimum can then easily be done using any complete tree search algorithm and evaluating the score at each leaf of the tree to find the minimum.

The first optimization to the brute-force algorithm is to compute the score during the exploration of the nodes in the tree, while doing a depth first search. If a node in the tree has a higher score than the minimum so far, the whole branch can be pruned. Further pruning is possible by using the monotonicity property (section~\ref{sec:monotonicity}). If the number of scenes discovered by a query is higher than in the previous step, the whole branch can also be pruned. The property also allows to determine a lower bound on the score of a branch, which allows pruning earlier when compared to the current minimum score.

Another optimization considers the possibility to permute the elements making up the queries as described in section~\ref{sec:uniqueness}. When the search tree is pruned by only allowing one of these permutations, the number of sequences to be checked is a factor $n!$ lower.

Additionally, at some point, all the remaining queries remove exactly one scene each. All possible permutations of these queries will lead to the same score. When looking for only one optimum the search can just take any random permutation of the remaining queries at this point or directly compute all permutations without recalculating the score for each of them.

Although the required time to find the optimal sequence drastically decreases, these optimizations do not influence the overall computational complexity, making it infeasible to compute the optimal sequence for $n>6$, $k>3$ on computers that were available to us at the time of writing.

\subsection{Goal-driven Algorithms}

In this section we propose two algorithms that allow us to compute sub-optimal solutions for higher values of ($n$, $k$). Both algorithms use heuristics to assign a score  $\delta(q_i)$ to a candidate query $q_i$ that indirectly estimates how well it contributes to minimize the sequence's overall time to discovery. Then, the query with the highest $\delta(q_i)$ is executed.

\subsubsection{Greedy Scene Elimination (GSE)}
\label{sec:GSE}

In this greedy algorithm, a sequence is built by simply selecting the queries $q_i$ that discovers the highest number of scenes $\delta_{\text{GSE}}(q_i) = D(q_i)$ at each step. This approach implicitly produces sequences having the monotonicity property described in section~\ref{sec:monotonicity}. Ties can be resolved in different ways. Our current implementation chooses the query that comes first in lexicographic order. We found that the choice influences the final score, which is evidence for the non-optimality of the algorithm. 
Initially, all the scenes are generated and stored. Then, with each query, the discovered scenes are eliminated from the list of scenes, hence the name of the algorithm. This equals following just one path to a single leaf in the tree search algorithm. The computational complexity of this algorithm is therefore lowered to $O\left(Nk |S|\right)$, while the memory consumption is $O\left(n |S|\right)$ to store the scenes that need to be checked.

While the branch and prune algorithm only allowed us to compute the optimal solution for $n < 7$, the greedy scene elimination algorithm allowed us to compute solutions for $n < 32$ within a few days.

\subsubsection{Minimally Intersecting Subsets (MIS)}

Given the computational requirements of GSE, a simpler scoring function is defined by considering the overlap between the candidate $q_i$ and each one of queries previously executed. For each $q_i$, we select the one with the highest:
$$\delta_{\text{MIS}}(q_i) = - \sum_{j=1}^{i-1}2^{|q_i \cap q_j|} - 1$$
The term $2^{|q_i \cap q_j|}$  corresponds to the number of subsets in the intersection $q_i \cap q_j$. We exclude the empty subset from the count to avoid undesirable results: two queries with $|q_i \cap q_j|=0$ would be equivalent to one with $|q_i \cap q_j|=1$. The same procedure as in GSE is used to break ties.

The score $\delta_{\text{MIS}}(q_i)$ can be seen as an approximation to $\delta_{\text{GSE}}(q_i)$. Let $d(q_i)$ be the set of scenes $\subseteq S$ that could be discovered by $q_i$, we can express  $\delta_{\text{GSE}}(q_i)$ as:
\begin{IEEEeqnarray*}{ll}
\hspace{2em}&\hspace{-2em} \delta_{\text{GSE}}(q_i) = D(q_i) 
= | d(q_i) |  \\
& - \sum_{j=1}^{i-1}  | d(q_i) \cap d(q_j) | \\
& +  \sum_{j=1}^{i-1}\sum_{k=1}^{j-1}  | d(q_i) \cap d(q_j) \cap d(q_k) | \\
& -  \sum_{j=1}^{i-1}\sum_{k=1}^{j-1} \sum_{l=1}^{k-1}  | d(q_i) \cap d(q_j) \cap d(q_k)\cap d(q_l) | \\
& +  \ldots \IEEEyesnumber \label{eq:sc_removed} \\ 
\end{IEEEeqnarray*}
where the first term is equal for all $q$ and thus removing it results in an equivalent score function. Here we consider an approximation to $\delta_{\text{GSE}}$  by only considering the second term: 
$$\delta_{\text{MIS}^*}(q_i) = - \sum_{j=1}^{i-1}  | d(q_i) \cap d(q_j) | $$
It can be shown that the size of the intersection between the scenes removed by $q_i, q_j$ can be expressed more easily as:
$$\delta_{\text{MIS}^*}(q_i) = -\sum_{j=1}^{i-1} 2 ^{n - 2k + |q_i \cap q_j|}$$
This corresponds to $\delta_{\text{MIS}}(q_i)$, except for a factor $2^{n - 2k}$ and the exclusion of the empty subset in the count as explained above. Increasingly better approximations to GSE could be obtained by considering more and more terms of $\eqref{eq:sc_removed}$, an option we have not studied in this paper.

\subsection{Stateless Sequence Generation Algorithms}

The algorithms described in this section do not explicitly consider the objective in equation~\eqref{eq:ttd} to generate $k$-subsets. Thus they do not require to store the previously executed queries or the remaining stars, not having a representation of the current state of the problem. This leads to implementations with lower computational and memory complexity. 

These algorithms can be seen as analogous to those used to solve the minimal change ordering problem such as the revolving door algorithm~\cite{cages}, the Eades-McKay algorithm~\cite{matters}, and Chase's sequence~\cite{knuth}. Trying to minimize the distance between consecutive $k$-subsets, they perform badly in our case, as close, similar subsets discover fewer scenes.

\subsubsection{Generalized Pattern Shifting Algorithm}

All the previous methods require the evaluation of a score function for every query, which requires a high computational cost even for simple functions as $n$ and $k$ increase. Those algorithms cannot run in real-time on a star tracker. In \cite{MortariEtAl2004}, Mortari et al. propose an algorithm for producing $k=3$ subset sequences. The algorithm is notable for the low computational complexity with which it generates sequences, while maintaining some degree of diversity among its sequentially produced subsets.
We present here a generalization of the pattern shifting algorithm to arbitrary $k$, and arbitrary reference sequences.

Let $E = \{e_0, \ldots, e_{n-1}\}$ represent the $n$ elements being grouped into $k$-element subsets.
A reference sequence generator is made to iteratively produce its sequence of $(k-1)$-element subsets from $E \setminus \{e_0\}$.
Each subset $q'$ produced this way is extended with $e_0$, so as to compose a valid $k$-element subset: $q_i = \{e_0\} \cup q'$.
Multiple variants of this subset are then produced, by incrementing/shifting subset elements: $q_{i+1} = \{e_{m+1} | e_m \in q_i\}$.
This will lead to a sub-sequence $\{q_i, \ldots, q_j\}$ of $Q$, ending at the subset $q_j$ for which $e_{n-1} \in q_j$. Should this condition hold already for $q_i$, the sub-sequence will then include a single subset.

In the following sections we take lexicographic as the reference sequence generator, which, for $k=3$, then perfectly replicates the sequences produced by the algorithm given in \cite{MortariEtAl2004}. However, the formulation presented above allows for other reference generators to be considered. It can even use itself recursively stopping at $k=0$ which simply results in a single empty subset.
Naturally, the quality of the sequences generated by this algorithm are closely tied to the quality of the sequences produced by the reference generator. Exchanging lexicographic for revolving door \cite{cages}, for instance, will improve the overall quality.

\subsubsection{Unranking Algorithms}

We here consider a class of sequence generators that rely on unranking functions to produce their sequences. These are functions that, for a given rank $r$, will generate the $(r+1)$-th subset in some reference sequence. In \cite{cages}, unranking functions are provided for the lexicographic, co-lexicographic and revolving door sequence generators.

By relying on an unranking function, the problem of generating a sequence is recast as the problem of how to best generate a sequence of $N$ ranks in $\{0, \ldots, N-1\}$. A pseudorandom number generator (PRNG) can be used to produce these ranks. If the uniqueness of ranks is enforced, for instance via a full cycle PRNGs such as linear feedback shift registers, we then obtain a pseudorandom sequence without repetitions, effectively a pseudorandom permutation of the $N$ subsets. In terms of scalability, unranking functions enable such sequences to be generated without the need to obtain ahead the full sequence $Q$ and then perform a shuffling operation in memory over it.

This approach can however be improved through analysis of the reference sequences for which the unranking functions are designed. Lexicographic, co-lexicographic and revolving door are all sequences where similar subsets can be found clustered with similar ranks. A logical approach is then to seek a process that generates ranks that are locally as dissimilar as possible, and preserves this property throughout the full sequence. Techniques to achieve that can be found in the domain of low-discrepancy, or quasi-random sequence generators. We introduce here an approach, under the name ``base unrank", that is based on the van der Corput sequence \cite{corput}.

Let $L = \lceil log_b{N} \rceil$ represent the number of digits required in order to count from 0 to $N-1$ in base $b$.
The sequence of ranks is produced by counting from $0$ to $b^L - 1$, in base $b$, while prioritizing increments to the \emph{most significant digit}. This is equivalent to a sequence that increments ranks as per the normal rules of addition, but at each step then reverses the $L$ digits. Ranks $r \geq N$ produced by this process are simply skipped in the sequence. This sequence generation equals using the van der Corput sequence from $0$ to $b^L - 1$ multiplied by $b^L$ as index for the unranking function.

Consider the example of base unranking for $n=5$, $k=3$ ($N=10$). Taking $b=2$ will result in a counting over $L=4$ digits. In binary, the sequence would be $\{0000, 1000, 0100, \ldots, 1111\}$, which would decode in base 10 to $\{0, 8, 4, 12, 2, 10, 6, 14, 1, 9, 5, 13, 3, 11, 7, 15\}$. Skipping the ranks $r \geq N$ would then lead to the sequence $\{0, 8, 4, 2, 6, 1, 9, 5, 3, 7\}$.
These ranks, produced sequentially, and mapped through the unranking function for some reference sequence such as the revolving door (the default used in the remainder of this paper), will generate a valid $k$-elements subset sequence with no repetitions, and little overhead (computationally, the costlier operation will be the unranking itself).

\iftrue

\subsubsection{Random Queries}

For completeness, we also consider the expected score of purely random queries, selecting $k$-subsets at uniform probability from the set of all possible $k$-subsets. This means, that $Q$ would contain repetitions and has infinite length. In practice it would be necessary to abort the discovery at some point to avoid running an endless loop, when a scene cannot be discovered.
With a derivation similar to that of Eq.(\ref{eq:sigma}), we obtain:
\begin{equation}
   \expect{T(Q)} = \frac{1}{|S|} \sum_{t=k}^{n} \left( |S_t| \cdot \frac{N}{\binom{t}{k}} \right)
\end{equation}
\fi

\section{Results}

\begin{figure}[t] 
    \centering
    \subfloat[][Optimal]{
        \includegraphics[width=.1\textwidth]{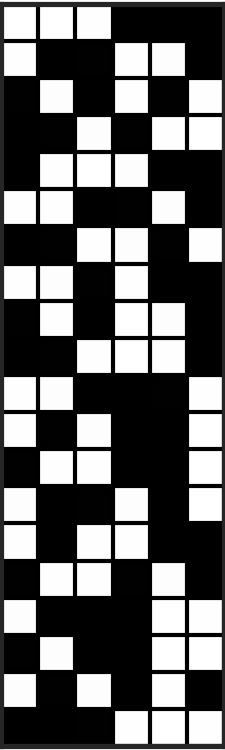}\vspace{1em}
    }
    \subfloat[][GSE]{
        \includegraphics[width=.1\textwidth]{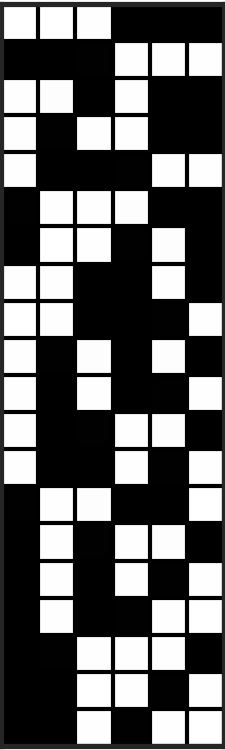}\vspace{1em}
    }
    \subfloat[][Scenes Discovered]{
        \includegraphics[width=.2\textwidth]{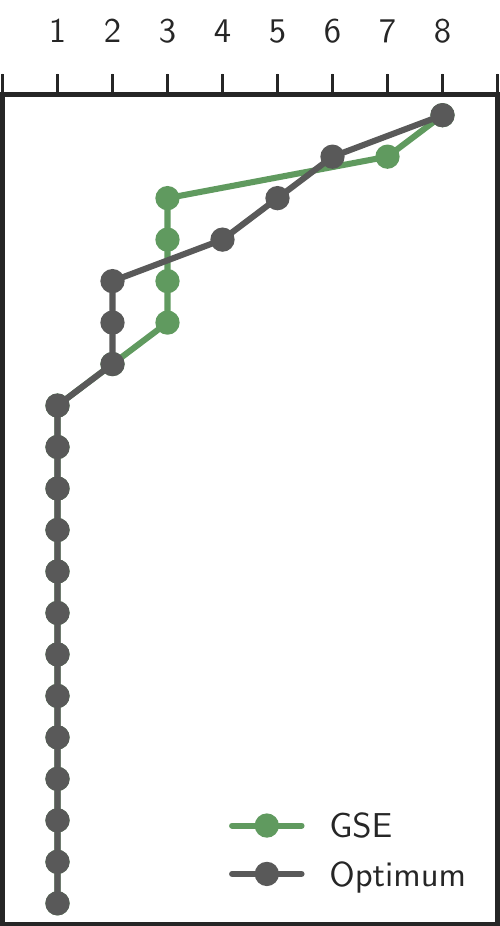}
    }
    
    \caption{For $n=6$ and $k=3$ an optimal solution (a) cannot be found by the greedy scene elimination (b) algorithm. This is because the second query of the optimal solution discovers less scenes than GSE at this step, but is then able to discover more scenes with the subsequent two queries.}
    \label{fig:sequences_6_3}
\end{figure}

\begin{figure*}[t]
    \centering
    \includegraphics[width=\textwidth]{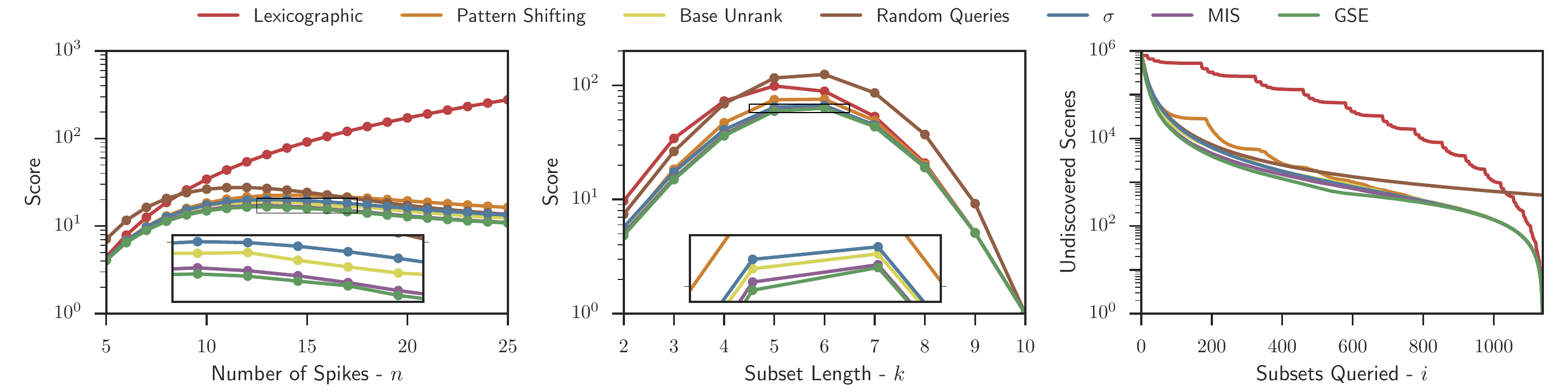}
    \caption{Comparison of the performance of the algorithms presented. The left and center panel show the score $T(Q)$ of the sequences $Q$ generated by the algorithms depending on the number of spikes $n$ (with $k=3$, left panel) and the subset length $k$ (with $n=10$, center panel). The right panel shows the number of scenes left to discover for the sequences $Q$ after every query $q_i$ for $i$ from 0 to $N-1$ for the scenario with $n=20$ and $k=3$.}
    \label{fig:results}
\end{figure*}

We test all the proposed algorithms in scenarios with different values of ($n$, $k$), reporting the average time to discovery $T(Q)$ of each case. We pay particular attention to the $k=3$ cases, due to their relevance to star trackers.
Additionally, we analyse the sequences produced in some of these cases to study how the different algorithms behave.

The branch and prune algorithm can find an optimum 
sequence for scenarios up to $n = 6$. In the case of $n=6, k=3$, it is interesting to note that there is an overlap between the first two queries,  $[0, 1, 2]$ followed by $[0, 3, 4]$ as illustrated in Figure~\ref{fig:sequences_6_3}. It is interesting to note that such sequence is not generated by GSE and MIS algorithms, given that they try to select the query discovering the maximum amount of scenes at each step, thus minimizing the overlap with respect to the previous sequence.

\begin{figure*}[p] 
\captionsetup[subfigure]{labelformat=empty} 
\captionsetup[subfigure]{justification=centering}
    \centering
    \subfloat[][Lexicographic]{
        \includegraphics[width=.099\textwidth]{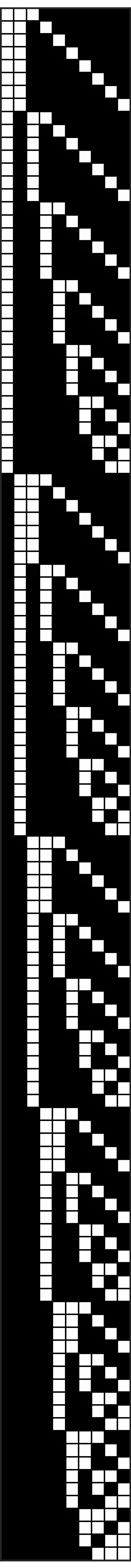}
    }
    \hspace{10pt}
    \subfloat[][Revolving Door]{
        \includegraphics[width=.099\textwidth]{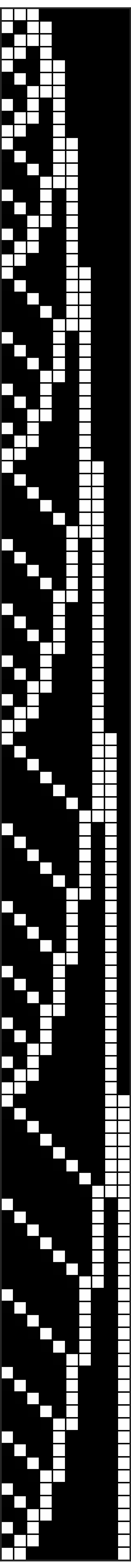}
    }
    \hspace{10pt}
    \subfloat[][Pattern Shifting]{
        \includegraphics[width=.099\textwidth]{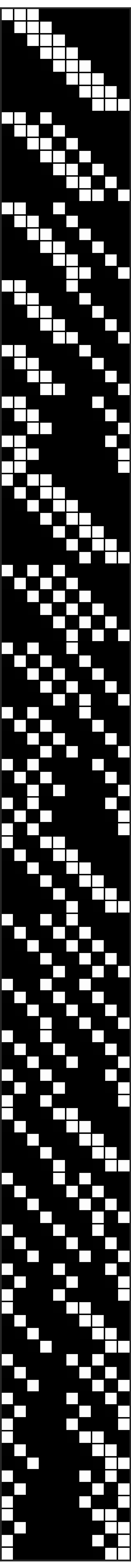}
    }
    \hspace{10pt}
    \subfloat[][Greedy Scene\\Elimination]{
        \includegraphics[width=.099\textwidth]{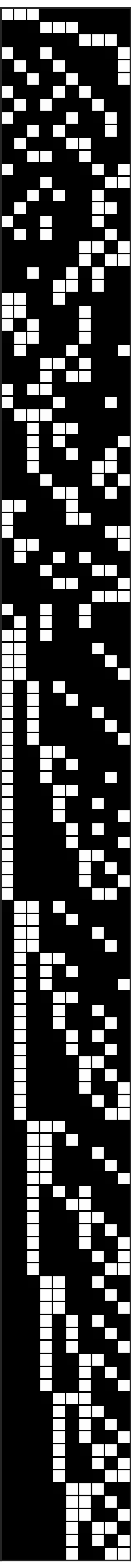}
    }
    \hspace{10pt}
    \subfloat[][Minimally\\Intersecting\\Subsets]{
        \includegraphics[width=.099\textwidth]{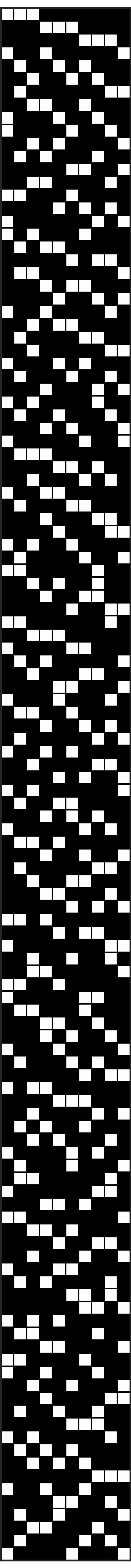}
    }
    \hspace{10pt}
    \subfloat[][Base Unrank]{
        \includegraphics[width=.099\textwidth]{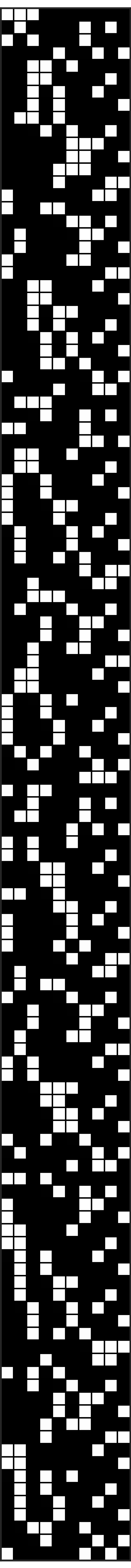}
    }
    \caption{Visualisation of the sequences generated by the algorithms implemented in this paper, for $n=10$ and $k=3$.}
    \label{fig:sequences_10_3}
\end{figure*}

As a reference, we include the sequences generated by the algorithms for $n=10, k=3$ in Figure  \ref{fig:sequences_10_3}. Note that for GSE, in this case, after 48 queries, all the sequences remove one scene, thus resulting in a lexicographic ordering of the remaining queries.

Figure~\ref{fig:results} shows the results of the algorithms in scenarios with different ($n$, $k$). Excluding the random queries case, the remaining algorithms can be ordered according to their performance from best to worst as follows:
$$\text{GSE} \prec \text{MIS}  \prec \text{base unrank}  \prec \text{pattern-shifting} \prec \text{lex}$$
It is interesting to note that the pattern-shifting and lexicographic (lex) algorithms exhibit a performance worse than $\sigma$. Although the pattern-shifting algorithm has a remarkably simple implementation and has low complexity both in terms of computations and memory, a random ordering is likely to perform better. Moreover, base-unrank offers even better results with comparable costs.

Another interesting point is that, for $n >> k$, as repetitions are less likely, the performance of the random queries converge to $\sigma$, as can be seen in Figure~\ref{fig:results} in the left panel. In this case, it becomes better than the pattern shifting algorithm for $n \ge 17$.

The magnified areas in Figure~\ref{fig:results} show that the margin between the algorithms is relatively small, especially between GSE and MIS. A similar difference is expected between GSE and the optimal solution.

To investigate the behavior of the algorithms, Figure~\ref{fig:results} (right panel) shows how many scenes are undiscovered in the same scenario ($n=20$, $k=3$), while iterating over the queries in the sequence. The lexicographic algorithm shows a repeated pattern where less and less scenes are removed up to some point, where suddenly more scenes are removed again. The pattern shifting algorithm shows a similar behaviour, though less pronounced. Interestingly, the lines of $\sigma$ and the base unrank algorithm are intersecting. The better score of the base unrank algorithm seems to be caused by it removing more scenes in the beginning.

Of all algorithms without repetitions of queries, we saw that the lexicographic sequence performs worst. We find that the lexicographic sequence seems to be the worst case sequence for the described optimization problem. For $n < 6$, the lexicographic sequence is always among the sequences that maximizes the score as we found during the brute-force searches. Reversing the choice of GSE to choose the query that removes the least amount of scenes also yields a sequence with exactly the same score of the lexicographic sequence and showing the exact same behavior in the plot shown in Figure~\ref{fig:results}.

The best solutions generated for each one of the scenarios considered in this paper have been made available online\footnote{\url{http://www.esa.int/gsp/ACT/ai/projects/star_trackers.html}}.

\section{Conclusion}

We consider a star identification problem and we map it to a $k$-subsets (queries) optimal ordering problem. We provide an in depth analysis of its structure proving interesting mathematical properties that are used in the design and assessment of solution algorithms. A number of algorithms with different complexity and performance, covering a wide spectrum of the parameters $n$ and $k$, is proposed and proved to advance the state-of-the-art in star identification research. For small $n<7$ we are able to provide the optimal solutions, while for higher $n$ ($< 32$) our algorithms are only able to compute and score suboptimal solutions. We released our best solutions online\footnotemark[1]. The problem complexity is such that for even higher $n$ it is extremely challenging to even compute the considered scoring function. Nevertheless, we presented a class of algorithms with polynomial complexity in the output size $N$, producing better sequences than the average random sequence at least for the parameter range that could be tested.





\bibliographystyle{IEEEtran}
\bibliography{star_tracker}

\end{document}